\documentclass{svjour3}

\usepackage{amscd,amsfonts,amssymb,amsmath,latexsym,array,hhline,xcolor,graphicx}

\usepackage{latexsym}
\usepackage{times}

\newcommand\F{\mbox{I\kern-2pt F}}
\newcommand{\RNumb}[1]{\uppercase\expandafter{\romannumeral #1\relax}}

\newcommand\beq{\begin{equation}}
\newcommand\eeq{\end{equation}}
\newcommand\bea{\begin{eqnarray}}
\newcommand\eea{\end{eqnarray}}
\newcommand\bean{\begin{eqnarray*}}
\newcommand\eean{\end{eqnarray*}}

\begin{document}

\title{On the Application of Laplace Transform to the Ruin Problem with Random Insurance Payments and Investments in a Risky Asset}

\author{ Viktor Antipov } 
\date{}
\institute{  \at Lomonosov Moscow State University and ``Vega" Institute, Moscow, Russia \\
\email{stayaptichek@gmail.com}
}


\titlerunning{On the Application of Laplace Transform to the Ruin Problem with Investments}

\maketitle

\begin{abstract}
This paper considers the ruin problem with random premiums, whose densities have rational Laplace transforms, and investments in a risky asset whose price follows a geometric Brownian motion. The asymptotic behavior of the ruin probability for large initial capital values is investigated. 
\end{abstract}



%\begin{document}

%\author{Antipov V.~A.\Addressmark[1]}
%
%\Addresstext[1]{Lomonosov Moscow State University, Faculty of Mechanics and Mathematics, Department of Probability Theory, Moscow, Russia; e-mail: vantipov@nes.ru. Vega Science Development Support Foundation.}
%
%
%\udk{}

%\dedicatory{To the sixtieth anniversary}

%\title{On the Application of Laplace Transform to the Ruin Problem with Random Insurance Premiums and Investments in a Risky Asset}

%\markboth{Antipov V.~A.}{On the Application of Laplace Transform to the Ruin Problem with Investments}

% DON'T FORGET TO RESTORE!
%\maketitle

%\begin{fulltext}

%\begin{abstract}
%This paper considers the ruin problem with random premiums, whose densities have fractional-rational Laplace transforms, and investments in a risky asset whose price follows a geometric Brownian motion. The asymptotic behavior of the ruin probability for large initial capital values is investigated. 
%\end{abstract}

% DON'T FORGET TO RESTORE!
%\begin{keywords}
%ruin probability, actuarial models with investment, integro-differential equations, Laplace transform.
%\end{keywords}
%%%%%%%%%%%%%%%%%%%%%%%%%%%%%%%%%%%%%%%%%%%%%%%%%%%%%
%%%%%%%%%%%%%%%%%%%%%%%%%%%%%%%%%%%%%%%%%%%%%%%%%%%%%

\section{Introduction}
\label{s1}
One of the key questions in actuarial risk theory is the problem of finding the ruin probability or, at least, determining its asymptotic behavior as the initial capital \( u \) tends to infinity. In the classical Cramér-Lundberg model, it is assumed that an insurance company holds its capital in a risk-free asset. In particular, for light-tailed jump distributions, the ruin probability decays exponentially as the initial capital grows. In later models, this assumption is generalized—the insurance company reinvests its capital (or a portion thereof) into a risky asset, whose price dynamics follow a specific law. It turns out that adding investments to the model leads to a qualitatively different behavior of the ruin probability as the initial capital tends to infinity: the ruin probability decays according to a power law; in cases where the volatility of the risky asset is sufficiently high, ruin occurs with probability one. 

One approach to studying the behavior of the ruin probability is based on the analysis of integro-differential equations \cite{KPergFr2002}--\cite{KPerg2022}, \cite{BelkKonKur2012}--\cite{BelkKonSlav2019}. The smoothness of the ruin probability is investigated in \cite{KPergFr2002}--\cite{KPukh2020}, \cite{AK2024}. Asymptotic analysis for exponentially distributed insurance payments is presented in \cite{KPergFr2002}--\cite{KPukh2020}. The use of Laplace transform techniques for a Cramér--Lundberg model with investments in a risky asset is demonstrated in \cite{AlbrConstThom2012}--\cite{ConstThom2005}. One of the main results for the model with investments in a risky asset, described by a geometric Brownian motion with drift \( a \) and volatility \( \sigma > 0 \), and exponentially distributed jumps, is the asymptotic behavior of the ruin probability as the initial capital tends to infinity: the ruin probability decays as \( Cu^{-\beta}, C > 0 \), where \( \beta = 2a/\sigma^2 - 1 > 0 \), and ruin occurs with probability one for \( \beta \leq 0 \).  

The works \cite{KPergFr2002}, \cite{KPerg2016}, \cite{AlbrConstThom2012}, \cite{ConstThom2005} consider models with investments in a risky asset and one-sided insurance payments. The paper \cite{KPukh2020} studies the asymptotic behavior of the ruin probability in a model with two-sided jumps having exponential distributions. 
This paper examines a model with investments in a risky asset, allowing for both positive and negative insurance payments with density functions having rational Laplace transforms. The asymptotic behavior of the ruin probability as the initial capital tends to infinity is derived for this model. 

The paper is structured as follows. Section \ref{s2} describes the model. Section \ref{s3} presents auxiliary results on integral operators arising in the integro-differential equations describing the ruin probability. Section \ref{s4} contains the main results on the asymptotic behavior of the ruin probability in the model with risky investments.
 
\section{Model Description}
\label{s2}
Let a probability space \( (\Omega, \mathcal{F}, \mathbf{P}) \) with filtration \( \mathbf{F} = (\mathcal{F}_t)_{t \geq 0} \) be given, along with independent Wiener process \( W = (W_t)_{t \geq 0} \) and compound Poisson process \( P = (P_t)_{t \geq 0} \) with drift \( c \in \mathbb{R}\backslash \{0\} \). Define the Poisson random measure \( p_{P}(dt, dx) \) of process \( P \) with compensator \( \Pi_{P}(dx)dt \), where \( \Pi_{P}(\mathbb{R}) < \infty \). Further, assume that \( \Pi_{P}(\mathbb{R}_{+}) > 0 \) and \( \Pi_{P}(\mathbb{R}_{-}) > 0 \). 

Following \cite{KPukh2020}, consider a model of an insurance company investing its reserve in risky assets, whose price is described by a geometric Brownian motion 
\[ S_t = S_0 e^{(a - \sigma^2/2)t + \sigma W_t}, \quad a \in \mathbb{R}, \quad \sigma > 0. \]
The dynamics of the company's reserve are described by the process \( X = X^{u} \) of the form 
\[ X_t = u + \int_{0}^{t} X_s dR_s + P_t, \qquad t \geq 0, \]
where \( R_t = \log(S_t/S_0) = a t + \sigma W_t \) is the process of relative change in the risky asset price, the process \( P \) is a compound Poisson process with drift \( c \in \mathbb{R}\backslash \{0\} \), reflecting the "insurance" activity of the company, and the value \( u > 0 \) represents the initial capital of the company. The capital process \( X \) can also be represented in "differential" form as \( dX_t = X_t dR_t + dP_t, X_0 = u \), where 
\[ dP_t = cdt + \int x p_{P}(dt, dx), \quad P_0 = 0. \]
In actuarial literature, the process \( P \) is typically represented more intuitively as 
\[ P_t = ct + \sum_{i=1}^{N_t} \xi_i, \]
where \( N_t := p_{P}([0,t] \times \mathbb{R}) \) is a Poisson process with intensity \( \lambda = \Pi_{P}(\mathbb{R}) \), and the random variables \( \xi_i \) are independent and identically distributed with distribution function \( F(dx) = \Pi_{P}(dx) / \lambda \). In our case, it will be more convenient to represent the process \( P \) as 
\[ P_t = ct + \sum_{i=1}^{N^2_t} \xi^{2}_i - \sum_{i=1}^{N^1_t} \xi^{1}_i, \]
where \( N^1, N^2 \) are Poisson processes with intensities \( \lambda_1 = \Pi_{P}(\mathbb{R}_{-}) \) and \( \lambda_2 = \Pi_{P}(\mathbb{R}_{+}) \), respectively, and the random variables \( \xi^1_i, \xi^2_i \) are positive with the distribution functions \( F_1(x) := \Pi_{P}((-x, 0]))/\lambda_1 \) and \( F_2(x) := \Pi_{P}((0, x]))/\lambda_2 \) for \( x > 0 \). 

Following \cite{AlbrConstThom2012}, assume that the laws of random variables \( \xi^1_i \) and \( \xi^2_i \) are absolutely continuous, and their density functions \( f_1 \) and \( f_2 \) satisfy the differential equations 
\begin{equation}
\label{densities_equation}
\mathcal{P}_1 \left(\frac{d}{dx} \right)f_1(x) = 0, \quad \mathcal{P}_2 \left(\frac{d}{dx} \right)f_2(x) = 0,
\end{equation}
where the differential operators \( \mathcal{P}_i, i = 1, 2 \), are given by 
\[ \mathcal{P}_{i} \left( \frac{d}{dx} \right) = \sum_{j=0}^{n} \alpha^{j}_{i} \dfrac{d^j}{dx^j}, \quad \alpha^j_i \in \mathbb{R}, \quad \alpha^{0}_i \neq 0, \]
with boundary conditions 
\begin{equation}
\label{bound_conds}
f_i^{(k)}(0) = f^k_i, \quad f^k_i \in \mathbb{R}, \quad k = 0, \dots, n-1.
\end{equation}
Also define the corresponding adjoint operators \( \mathcal{P}^{*}_i, i = 1, 2 \),  
\[ \mathcal{P}^{*}_{i} \left( \frac{d}{dx} \right) = \sum_{j=0}^{n} (-1)^j \alpha^{j}_{i} \dfrac{d^j}{dx^j}. \]

Let \( \widehat f \) denote the one-sided Laplace transform of function \( f \), defined as 
\[ \widehat f(s) = \int_{0}^{\infty} e^{-sx} f(x) dx. \]
Sometimes it will be more convenient to use the one-sided Laplace--Stieltjes transform 
\[ \widehat U(s) = \int_{0}^{\infty} e^{-sx} dU(x), \]
where \( U(dx) \) is a measure defined on \( \mathbb{R}_{+} \). 

Note that the density functions \( f_1(x) \) and \( f_2(x) \), satisfying (\ref{densities_equation}), have Laplace transforms representable as rational fractions \( p(s) / q(s) \), where \( p(s) \) and \( q(s) \) are finite-degree polynomials (see \cite{AsmAlbr}, Ch. I).

Introduce standard notation: \( \tau := \inf \{t \geq 0: X^{u}_t < 0\} \) (time of ruin), \( \Psi(u) := \mathbb{P}(\tau < \infty) \) (ruin probability), \( \beta := 2a/\sigma^2 - 1 \).

\section{Integro-Differential Equation}
\label{s3}
According to Proposition 6.1 in \cite{KPukh2020}, the ruin probability \( \Psi \) satisfies the integro-differential equation 
\begin{eqnarray}
\label{mainIDE}
\mathcal{L}(\Psi)(u) + \mathcal{I}(\Psi)(u) = 0,
\end{eqnarray}
where 
\[ \mathcal{L}(\Psi)(u) = \frac{1}{2} \sigma^2 u^2 \Psi''(u) + (au + c) \Psi'(u) - (\lambda_1 + \lambda_2) \Psi(u), \]
\[ \mathcal{I}(\Psi)(u) = \int \Psi(u + y) \Pi_{P}(dy). \]
It will be convenient to work with the following representation of the integral operator:
\[ \mathcal{I}(\Psi)(u) = \lambda_1 \mathcal{I}_1(\Psi)(u) + \lambda_2 \mathcal{I}_2(\Psi)(u) = \lambda_1 \int \Psi(u - y) f_1(y)dy + \lambda_2 \int \Psi(u + y) f_2(y) dy. \]

Assume further that for some \( \beta' \in (0, \beta \wedge 1) \) the inequality \( \mathbb{E} (\xi_1^{1})^{\beta'} < 1 \) holds. Then, the ruin probability satisfies the regularity condition 
\begin{equation}
\label{regularity}
\lim_{u \to \infty} \Psi(u) = 0
\end{equation}
(see Lemma 2.1 in \cite{KPukh2020}).

\begin{proposition}
\label{propos_1}
For \( \mathcal{I}_1(\Psi)(u) \), the following equation holds:
\[ \mathcal{P}_1\left(\frac{d}{du}\right) \mathcal{I}_1(\Psi)(u) = \sum_{k=0}^{n-1} \Psi^{(k)}(u) \sum_{i=k}^{n-1} \alpha^{i+1}_1 f^{i-k}_1. \]
\end{proposition}
\proof
It is easy to see that 
\[ \int \Psi(u - y) f'_1(y)dy = \frac{d}{du} \mathcal{I}_1(\Psi)(u) - f_1(0) \Psi(u), \]
\[ \int \Psi(u - y) f''_1(y)dy = \frac{d^2}{du^2} \mathcal{I}_1(\Psi)(u) - f_1(0) \Psi'(u) - f'_1(0) \Psi(u). \]
By induction, we obtain 
\begin{equation}
\label{rel_1}
\int \Psi(u - y) f^{(k)}_1(y)dy = \frac{d^k}{du^k} \mathcal{I}_1(\Psi)(u) - \sum_{i=1}^{k} f^{(i-1)}_1(0) \Psi^{(k-i)}(u), \quad 0 \leq k \leq n.
\end{equation}
Multiplying both sides of (\ref{rel_1}) by \( \alpha^k_1 \), accounting for the boundary conditions, and summing the terms, 
\[ \int \Psi(u - y) \mathcal{P}_1\left(\frac{d}{du}\right) f_1(y)dy = \mathcal{P}_1\left(\frac{d}{du}\right) \mathcal{I}_1(\Psi)(u) - \sum_{k=1}^{n} \alpha^{k}_1 \sum_{i=1}^{k} f^{i-1}_1 \Psi^{(k-i)}(u). \]
Due to (\ref{densities_equation}), the left-hand side of the last equality is zero. After grouping terms on the right-hand side, we obtain the required relation.  
\endproof

\begin{proposition}
\label{propos_2}
For \( \mathcal{I}_2(\Psi)(u) \), the following relation holds:
\[ \mathcal{P}^*_2\left(\frac{d}{du}\right) \mathcal{I}_2(\Psi)(u) = \sum_{k=0}^{n-1} (-1)^{k} \Psi^{(k)}(u) \sum_{i=k}^{n-1} \alpha^{i+1}_2 f^{i-k}_2. \]
\proof
The proof follows the same steps as the proof of Proposition \ref{propos_1}, except that relation (\ref{rel_1}) is rewritten as 
\[ \int \Psi(u + y) f^{(k)}_2(y)dy = (-1)^k \frac{d^k}{du^k} \mathcal{I}_2(\Psi)(u) - \sum_{i=1}^{k} (-1)^{k-i} f^{(i-1)}_2(0) \Psi^{(k-i)}(u). \]
\endproof
\end{proposition}

\section{Asymptotic Analysis of the Ruin Probability}
\label{s4}

The ruin probability \( \Psi(u) \) satisfies (\ref{mainIDE}), which we rewrite as 
\[ \frac{1}{2} \sigma^2 u^2 \Psi''(u) + (au + c) \Psi'(u) - (\lambda_1 + \lambda_2) \Psi(u) + \lambda_1 \mathcal{I}_1(\Psi)(u) + \lambda_2 \mathcal{I}_2(\Psi)(u) = 0. \]
Define the differential operator 
\[ \mathcal{T} = \mathcal{P}_1 \mathcal{P}^{*}_2 = \sum_{j, k = 0}^{n} (-1)^{k} \alpha_1^{j} \alpha_2^{k} \dfrac{d^{j + k}}{dx^{j + k}}. \]
Applying it to the equation above, we obtain 
\begin{equation}
\label{with_inv_high_ord_eq}
q_{2n+2}(u) \Psi^{(2n + 2)}(u) + q_{2n+1}(u) \Psi^{(2n+1)}(u) + \ldots + q_1(u) \Psi'(u) = 0,
\end{equation}
where  
\begin{eqnarray*}
q_j(u) &=& a_j u^2 + b_j u + c_j, \qquad 0 \leq j \leq 2n + 2, \\
a_j &=& \frac{\sigma^2}{2} \underset{k+m=j-2}{\sum_{k, m=0}^{n}} \alpha_1^{k} \alpha_2^{m} (-1)^m, \\
b_j &=& \underset{k+m=j-1}{\sum_{k, m=0}^{n}} (\sigma^2(k + m) + a) \alpha_1^{k} \alpha_2^{m} (-1)^m, \\
c_j &=& d_j + g_j + \underset{k+m=j}{\sum_{k, m=0}^{n}} (k + m) \left(\frac{\sigma^2}{2} (k + m) + (a - \frac{\sigma^2}{2}) \right)\alpha_1^{k} \alpha_2^{m} (-1)^m, \\
d_{j} &=& c \underset{k+m=j-1}{\sum_{k, m=0}^{n}} \alpha_1^{k} \alpha_2^{m} (-1)^m - (\lambda_1 + \lambda_2) \underset{k+m=j}{\sum_{k, m=0}^{n}} \alpha_1^{k} \alpha_2^{m} (-1)^m, \\
g_{j} &=& \lambda_1 \underset{k+m=j}{\sum_{k, m=0}^{n-1,n}} \alpha_2^m (-1)^m \sum_{i=k}^{n-1} \alpha_1^{i+1} f_1^{i-k} 
%\\
%&+& 
+\lambda_2 \underset{k+m=j}{\sum_{k, m=0}^{n-1,n}} \alpha_1^m (-1)^k \sum_{i=k}^{n-1} \alpha_2^{i+1} f_2^{i-k}. 
\end{eqnarray*}
In the formulas above and hereafter, we assume \( \sum_{\emptyset} = 0 \). 

Note that 
\[ a_{2n+2} = \dfrac{\sigma^2}{2} \alpha_1^{n} \alpha_2^{n} (-1)^n > 0, \quad b_{2n+2} = c_{2n+2} = 0, \quad c_{2n + 1} = c \alpha_1^n \alpha_2^n (-1)^n \neq 0. \]
The function \( u q_{2n+1}(u) / q_{2n+2}(u) \) is not holomorphic at \( u = \infty \). Therefore, \( u = \infty \) is an irregular singular point of equation (\ref{with_inv_high_ord_eq}), and the Frobenius method \cite[Ch. 1, §2]{Fedoryuk} is not applicable. Introducing the notation \( G = \Psi' \) and applying the Laplace transform to the resulting equation, using the properties 
\begin{eqnarray*}
\widehat{G^{(k)}}(s) &=& s^k \widehat G(s) - \sum_{i=1}^{k} s^{k-i} G^{(i-1)}(0^{+}), \\
\widehat{u G^{(k)}}(s) &=& - s^k \widehat{G}'(s) - k s^{k-1} \widehat{G}(s) + \sum_{i=1}^{k-1} (k-i) s^{k-i-1} G^{(i-1)}(0^{+}), \\
\widehat{u^2 G^{(k)}}(s) &=& s^k \widehat{G}''(s) + 2k s^{k-1} \widehat{G}'(s) + k(k-1) s^{k-2} \widehat{G}(s) \\
&& - \sum_{i=1}^{k-2} (k-i)(k-i-1) s^{k-i-2} G^{(i-1)}(0^{+}),
\end{eqnarray*}
we arrive at the equation 
\begin{eqnarray}
\label{inv_lapl_tr_eq}
p_{2n+1}(s) \widehat{G}''(s) + l_{2n+1}(s) \widehat{G}'(s) + r_{2n+1}(s) \widehat{G}(s) = v_{2n}(s),
\end{eqnarray}
where 
\begin{eqnarray}
\label{first_coef}
p_{2n+1}(s) &=& \sum_{i=1}^{2n+1} a_{i + 1} s^{i}, \qquad v_{2n}(s) = \sum_{i=0}^{2n} \widetilde v_{i+1} s^{i}, \\
l_{2n+1}(s) &=& \sum_{i=0}^{2n+1} \widetilde b_{i+1} s^i, \qquad \widetilde b_{i} = \begin{cases} 2i a_{i+1} - b_i, \quad 1 \leq i \leq 2n+1, \\ -b_{2n+2}, \qquad i = 2n + 2, \end{cases} \\
\label{last_coef}
r_{2n+1}(s) &=& \sum_{i=0}^{2n+1} \widetilde c_{i+1} s^i, \qquad \widetilde c_{i} = \begin{cases} i(i+1)a_{i+2} - 2i b_i + c_i, \quad 1 \leq i \leq 2n, \\ -2(2n+1)b_{2n+2}, \qquad i = 2n + 1, \\ c_{2n+1}, \qquad i = 2n + 2.
\end{cases}
\end{eqnarray}

Note that explicit formulas for the coefficients \( \widetilde v_{i} \) are quite cumbersome and omitted, as their values do not affect the subsequent reasoning. Thus, the following proposition is proven. 
\begin{proposition}
The Laplace transform \( \widehat G \) of the function \( \Psi' \) satisfies the differential equation (\ref{inv_lapl_tr_eq}), whose coefficients are given by (\ref{first_coef}) - (\ref{last_coef}). 
\end{proposition}

Assume \( 2a/\sigma^2 < 2 \). Note that \( s = 0 \) is a regular singular point of (\ref{inv_lapl_tr_eq}), since 
\[ \widetilde l_0 := \lim_{s \to 0} s \frac{l_{2n+1}(s)}{p_{2n+1}(s)} = \frac{\widetilde b_{1}}{a_2} = \frac{2a_2 - b_1}{a_2} = 2 - \frac{b_1}{a_2} = 2 - \frac{2a}{\sigma^2}, \]
\[ \widetilde r_0 := \lim_{s \to 0} s^2 \frac{r_{2n+1}(s)}{p_{2n+1}(s)} = \lim_{s \to 0} s \frac{\widetilde c_1}{a_2} = 0. \]

Apply the Frobenius method \cite[Ch. 1, §2]{Fedoryuk} to solve (\ref{inv_lapl_tr_eq}). The indicial equation is 
\[ \rho (\rho - 1) + \widetilde l_0 \rho + \widetilde r_0 = 0. \]
It is easy to see that the solutions of the indicial equation are 
\[ \rho_1 = 0, \qquad \rho_2 = \frac{2a}{\sigma^2} - 1. \]
Then, for \( |s| < \delta \), the fundamental system of solutions to (\ref{inv_lapl_tr_eq}) is given by 
\begin{equation}
\label{fss}
\widehat G_1(s) = s^{\rho_1} \gamma_1(s), \qquad \widehat G_2(s) = s^{\rho_2} \gamma_2(s),
\end{equation}
where \( \gamma_i(s) \) are holomorphic functions for \( |s| < \delta \) and \( \gamma_i(0) \neq 0 \), \( i = 1, 2 \). 
The general solution of the homogeneous equation is 
\[ \widehat G^{\text{oo}}(s) = c_1 \widehat G_1(s) + c_2 \widehat G_1(s), \qquad c_1, c_2 \in \mathbb{R}. \]
Using the method of variation of constants, a particular solution is 
\[ \widehat G^{\text{ps}}(s) = - \widehat G_1(s) \int_{0}^{s} \dfrac{\widehat G_2(t) v_{2n}(t)}{p_{2n+1}(t) W(t)} dt + \widehat G_2(s) \int_{0}^{s} \dfrac{\widehat G_1(t) v_{2n}(t)}{p_{2n+1}(t) W(t)} dt, \]
where \( W(s) \) is the Wronskian of the fundamental system of solutions (\ref{fss}). 
\[ W(s) = \widehat G_1(s) \widehat G'_2(s) - \widehat G_2(s) \widehat G'_1(s) = s^{\rho_1 + \rho_2 - 1} \widehat \gamma(s), \]
where \( \widehat \gamma(s) \) is a function holomorphic near 0. 

Due to the estimates 
\[ \frac{\widehat G_2(s)}{p_{2n+1}(s) W(s)} \sim \tilde c_1 \frac{s^{\rho_2}}{s^{\rho_1 + \rho_2}} = \tilde c_1, \qquad \frac{\widehat G_1(s)}{p_{2n+1}(s) W(s)} \sim \tilde c_2 \frac{s^{\rho_1}}{s^{\rho_1 + \rho_2 }} = \tilde c_2 s^{- \rho_2}, \quad s \to 0, \]
where \( \tilde c_1, \tilde c_2 \in \mathbb{R} \), it follows that 
\[ \widehat G^{\text{ps}}(s) = O(s), \quad s \to 0. \]
Note that due to the regularity condition (\ref{regularity}), we have \( c_1 = 0 \). Thus, 
\[ \widehat G(s) \sim \widetilde C s^{\rho_2}, \quad s \to 0, \]
and by the Tauberian theorem \cite[Ch. XIII, Theorem 4]{Feller}, 
\[ \Psi(u) \sim C u^{-\rho_2}, \quad u \to \infty, \]
where \( C > 0 \).
\begin{theorem}
Let \( \beta = 2a/\sigma^2 - 1 \), where \( \beta \in (0, 1) \). Assume the density functions \( f_i, i = 1, 2 \), satisfy (\ref{densities_equation}) with boundary conditions (\ref{bound_conds}) and there is \( \beta' \in (0, \beta \wedge 1) \) such that \( \mathbb{E} (\xi_1^{1})^{\beta'} < 1 \). Then \( \Psi(u) \sim C u^{-\beta} \), \( u \to \infty \) for some \( C > 0 \). 
\end{theorem}

%\end{fulltext}


\begin{thebibliography}{100}

\bibitem{KPergFr2002}
 A.~Frolova, Yu.~Kabanov, S.~Pergamenshchikov: 
In the insurance business risky investments are dangerous
Finance and Stochastics, {\bf 6}
227--235 (2002)


\bibitem{KPerg2016}
%\Hrefs{https://doi.org/10.1007/s00780-016-0292-4}
Yu.~Kabanov, S.~Pergamenshchikov:
In the insurance business risky investments are dangerous: the case of negative risk sums
Finance and Stochastics, 
{\bf 20},  2
355--379 (2016)
%\crossref{https://doi.org/10.1007/s00780-016-0292-4}

\bibitem{KPukh2020}
%\Hrefs{doi:10.1017/jpr.2021.74}
Yu.~Kabanov, N.~Pukhlyakov:
Ruin probabilities with investments: smoothness, IDE and ODE, asymptotic behavior. 
J. Appl. Probab.
{\bf 59}, 2
556--570 (2020)
%\crossref{doi:10.1017/jpr.2021.74}

\bibitem{AlbrConstThom2012}
%\Hrefs{https://api.semanticscholar.org/CorpusID:10356190}
H.~Albrecher, C.~Constantinescu, E.~Thomann: 
Asymptotic results for renewal risk models with risky investments
Stochastic Processes and their Applications
{\bf 122}
3767--3789 (2012)
%\crossref{https://api.semanticscholar.org/CorpusID:10356190}

\bibitem{ConstThom2005}
C.~Constantinescu, E.~Thomann: 
Analysis of the ruin probability using Laplace transforms and Karamata Tauberian theorem. 
ARCH, Proceedings 39th Actuarial Research Conference, Iowa City, Iowa, 2004
{\bf 1} (2005)

\bibitem{AsmAlbr}
S.~Asmussen, H.~Albrecher:  
Ruin Probabilities. 
World scientific (2010)

\bibitem{KPerg2022}
%\Hrefs{https://doi.org/10.1007/s00780-022-00483-w}
Y.~Kabanov, S.~Pergamenshchikov:
On ruin probabilities with investments in a risky asset with a regime-switching price.
Finance and Stochastics,
{\bf 26},  877--897 (2022)
%\crossref{https://doi.org/10.1007/s00780-022-00483-w}

\bibitem{AK2024}
%\Hrefs{https://doi.org/10.3390/math12111705}
V.~Antipov, Y.~Kabanov: 
Ruin Probabilities with Investments in Random Environment: Smoothness. 
Mathematics, 
{\bf 12}, 11,  1705 (2024)
%\crossref{https://doi.org/10.3390/math12111705}

\bibitem{BelkKonKur2012}
%\Hrefs{https://doi.org/10.1134/S0965542512100077}
T.~Belkina, N.~Konyukhova, S.~Kurochkin:
%\yr 2012
Singular boundary value problem for the integrodifferential equation in an insurance model with stochastic premiums: Analysis and numerical solution
Comput. Math. and Math. Phys.,
{bf 52}, 1384--1416 (2012)
%\crossref{https://doi.org/10.1134/S0965542512100077}

\bibitem{BelkKonSlav2019}
%\Hrefs{https://doi.org/10.1134/S0965542519110022}
T.~Belkina, N.~Konyukhova, B.~Slavko: 
%\yr 2019
Solvency of an Insurance Company in a Dual Risk Model with Investment: Analysis and Numerical Study of Singular Boundary Value Problems. 
Comput. Math. and Math. Phys.
{\bf 59}, 
1904--1927 (2019)
%\crossref{https://doi.org/10.1134/S0965542519110022}

\bibitem{Feller}
W.~Feller:
%\yr 1971
An Introduction to Probability Theory and Its Applications, vol 2.
John Wiley and Sons Inc., New York (1971)

\bibitem{Fedoryuk}
M.~Fedoryuk: 
%\yr 1993
 Asymptotic Analysis: Linear Ordinary Differential Equations. Springer (1993)

\end{thebibliography}
\end{document}